\documentclass[a4paper,11pt,leqno]{article} 
\usepackage{amsmath,amscd,amsthm,amssymb,latexsym,oldgerm,eufrak,epsfig,graphics,array,enumerate} 
 
\addtocounter{equation}{0} 
\numberwithin{equation}{section}

\newcommand{\du}[2]{\mbox{\raisebox{0.5ex}{$#1$}}  
\mbox{\hspace*{-0.25em}} \diagup_{\displaystyle \mbox{\hspace*{-0.3em}} #2}}

\begin{document} 
 
\theoremstyle{theorem} 
\newtheorem{theo}{Theorem}
\theoremstyle{theorem} 
\newtheorem*{satz}{Satz} 
\theoremstyle{remark} 
\newtheorem*{rem}{Remark} 
\theoremstyle{remark} 
\newtheorem*{bew}{Beweis} 
\theoremstyle{corollary} 
\newtheorem*{coro}{Korollar} 
\theoremstyle{lemma} 
\newtheorem*{lem}{Lemma}
\theoremstyle{definition} 
\newtheorem*{defi}{Definition}
\theoremstyle{proposition} 
\newtheorem{prop}{Proposition}[section]
\theoremstyle{proposition} 
\newtheorem{prop1}{Proposition}[section]
\theoremstyle{theorem} 
\newtheorem{theo1}[prop1]{Theorem}
 
\title{ A Bernstein Theorem for Special Lagrangian Graphs} 
\author{by \\[2ex]
J\"urgen Jost \\
Max-Planck-Institute for Mathematics in the Sciences \\
Inselstrasse 22 - 26, D-04103 Leipzig \\[2ex]
and \\[2ex]
Y. L. Xin \\
Institute of Mathematics \\
Fudan University, Shanghai 200433, P.R. China
} 
\maketitle

\pagenumbering{arabic} 
\bigskip
\footnotetext{The second-named author is grateful to the
Max Planck Institute for Mathematics in the Sciences 
in Leipzig for its hospitality and support. }
\bigskip
\bigskip
\begin{abstract}
We obtain a Bernstein theorem for
special Lagrangian graphs in $ \mathbb C^n = \mathbb R^{2n} $ for
arbitrary $ n $ only assuming bounded slope but no quantitative restriction. 
\end{abstract} 
\bigskip
\bigskip
\section{Introduction}

Let $M$ be the graph in $ \mathbb C^n \cong \mathbb R^{2n} $ of a smooth map 
$ f : \Omega \to \mathbb R^n $,
with $ \Omega \subseteq \mathbb R^n $ an open domain. $ M $ is a Lagrangian 
submani\-fold of $ \mathbb C^n $ if and only if the matrix
$ \left( \frac{\partial f^i}{\partial x^j} \right)  $ is
symmetric. In particular, in that case if $ \Omega $ is simply connected, then
there exists a function $ F : \Omega \to \mathbb R $ with
$$
 \nabla F = f .
$$ 

A Lagrangian submanifold of $ \mathbb C^n $ is called special if it
is a minimal submanifold at the same time. 
In the above situation, the graph of $ \nabla F $ is a special
Lagrangian submanifold of $ \mathbb C^n $ if and only if for some constant $ \vartheta $,
\begin{equation}
\label{eq1.1}
 \mbox{Im } (\det (e^{i \vartheta} (I + i \mbox{ Hess } (F))) = 0
\end{equation}
(Im = imaginary part, $I$ = identity matrix, Hess 
$ F = ( \frac{\partial^2 F}{\partial x^i \partial x^j  })$). \\

Special Lagrangian calibration constitute an example of a calibrated
geo\-metry in the sense of Harvey and Lawson \cite{HL}.  
More recently, Strominger-Yau-Zaslow \cite{SYZ} established a 
conjectural relation of fibrations by special Lagrangian tori
with mirror symmetry. However, in general, some of these tori are
singular. More generally, understanding such fibrations systematically
means understanding the moduli space of special Lagrangian tori,
and for that purpose, one needs to study the possible singularities
of special Lagrangian submanifolds in flat space. By 
asymptotic expansions at singularities, the so-called blow-ups,
the study of singular special Lagrangian spaces is reduced to the study
of special Lagrangian cones. \\

In the theory of minimal submanifolds, there exists a close
link between the rigidity of minimal cones and Bernstein type theorems, 
saying that - under suitable boundedness
assumptions - entire minimal graphs are necessarily planar. 
Of course, the known Bernstein type theorems for
entire minimal graphs, in particular \cite{HJW} and \cite{JX},
apply here. Those results seem to be close to optimal 
already. It turns out, however, that under the Lagrangian 
condition, one may prove still stronger such results. This
is the content of the present paper. \\

Returning to (\ref{eq1.1}), the Bernstein question then is whether, 
or more precisely, under which conditions, an entire
solution has to be a quadratic polynomial. 

Fu \cite{Fu} showed that for $ n = 2 $, any solution
defined on all of $ \mathbb R^2 $ is harmonic or a quadratic polynomial. 
Before stating our results on this question, however, let
us briefly observe that the equation (\ref{eq1.1}) is similar 
to the Monge-Amp\`ere equation
\begin{equation}
\label{eq1.2}
 \det \left( \frac{\partial^2 F}{\partial x^i \partial x^j} \right) = 1
\end{equation}
that naturally arises in affine differential geometry. There,
one is interested in convex solutions, and Calabi \cite{C} 
showed that, for $ n \leq 5 $, any convex solution of (\ref{eq1.2}) that
is defined on all of $ \mathbb R^n $ has to be a quadratic
polynomial. Pogorolev subsequently extended this result to all $ n $
(\cite{P}, \cite{TW}).

In that direction, we have 
\begin{theo}
\label{theo1.1}
Let $ F : \mathbb R^n \to \mathbb R $ be a smooth function defined
on the whole $ \mathbb R^n $. Assume that the graph of $ \nabla F $ is a special Lagrangian 
submanifold $ M $ in $ \mathbb C^n = \mathbb R^n \times \mathbb R^n $, namely $ F $ satisfies 
equation (\ref{eq1.1}). If 
\begin{description}
\item{ (i)}	$ F $ is convex; 
\item{(ii)}	there is a constant $ \beta < \infty $ such that 
	\begin{align}
          \label{eq1.3}
		\Delta_F & \leq \beta , 
	\intertext{where}
          \label{eq1.4}
		\Delta_F & = \{ \det (I + ( \mbox{Hess }(F))^2 \}^{\frac{1}{2}} ,
	\end{align}
		then $F$ is a quadratic polynomial and $ M $ is an affine $n$-plane. 
\end{description}
\end{theo}
For the proof of this theorem, we shall use the same strategy as in our 
previous paper \cite{JX}, dealing with minimal graphs in general. \\

As $ M $ is a minimal submanifold of Euclidean space, 
by the theorem of Ruh-Vilms \cite{RV}, its Gauss map
is harmonic. The Gauss map takes its values in the 
Grassmannian $ G_{n, n} $ of $ n $-planes in $ 2n $-space. In order
to show that $ M $ is affine linear, we need to show
that the Gauss map is constant. 
The strategy of 
Hildebrandt-Jost-Widman \cite{HJW} then was to show that the 
image of the Gauss map is contained in some
geodesically convex ball, and to show a Liouville
type theorem to the extent that any such harmonic map
with values in such a ball is constant. The method
works optimally if we look at harmonic maps with values in
a space of constant sectional curvature, i.e. a sphere. 
In the case of higher codimension $ k $, the Grassmannian $ G_{n, k} $, 
however, does not have the same sectional curvature
in all directions anymore. The stra\-te\-gy of Jost-Xin \cite{JX}
then was to exploit the Grassmannian geometry more carefully
and to construct other geodesically convex sets
for which such a Liouville type theorem for 
harmonic maps still holds. This led to a considerable
strengthening of the Bernstein type theorems for minimal
graphs of higher codimension. Still, however, those
results do not yet imply the preceding theorem. 
We need to exploit the fact that $ M $ is not only minimal,
but also Lagrangian. In other words, its Gauss map takes
its values in a certain subspace of $ G_{n,n} $, namely the Lagrangian
Grassmannian $ L G_n $ of Lagrangian linear
subspaces of $ \mathbb R^{2n} $. $ LG_n $ is a totally geodesic 
subspace of $ G_{n,n} $, and so the Gauss map of $ M $ as a map into
$ LG_n $ is still harmonic. We can now exploit
the geometry of $ LG_n$ to construct suitable geodesically
convex subsets in that space and deduce a corresponding
Liouville type theorem. In that way, we shall 
show that the Gauss map of $ M $ is constant under the
conditions stated in Theorem \ref{theo1.1}, and so $ M $ is planar. \\

For proving that $ M $ is flat another possible approach is to study
its tangent cone $C \tilde M$ at infinity. In our situation
$ C \tilde M $ is a special Lagrangian cone. Its link is
a compact minimal Legendrian submanifold in $ S^{2n-1} $.
Thus, we prove Theorem \ref{theo1.2}. It is interesting in its
own right. 
\begin{theo}
\label{theo1.2}
Let $M $ be a simple (in the sense of \cite{HJW}) or compact
minimal Legendrian submanifold in $ S^{2n-1} $. Suppose
that there are a fixed $n$-plane $P_0 $ and some
$ \delta > 0 $, such that 
\begin{equation}
\label{eq1.5}
\langle P, P_0 \rangle \geq \delta
\end{equation}
holds for all normal $ n $-planes $ P $ of $ M $ in $ S^{2n-1} $. Then
$ M $ is contained in a totally geodesic subsphere of $ S^{2n-1} $.
\end{theo}
Although in general, any minimal $ 2$-sphere in a sphere
$ S^m $ is totally geodesic (see \cite{B}), there exist higher dimensional minimal submanifolds of $ S^m $
that are not totally geodesic and by \cite{LO}, we know that the analogue of Theorem \ref{theo1.2} does 
not hold for minimal submanifolds of spheres for arbitrarily
small values of $\delta$.
In fact, there even exist nontrivial
minimal Legendrian $ S^3 $'s in $ S^7 $ \cite{CDVV}, and already for
$ n = 2 $, one finds infinitely many different minimal Legendrian
two-dimensional Legendrian tori in $ S^5 $ \cite{H}. Thus, 
condition (\ref{eq1.5}) cannot be dropped even in the Legendrian case. \\

By using Theorem \ref{theo1.2} we can remove the convexity condition for the 
function $ F $ in Theorem \ref{theo1.1}. This is in fact our main result: 
\begin{theo}
\label{theo1.3}
Let $ F : \mathbb R^n \to \mathbb R  $ be a smooth function on the
whole $ \mathbb R^n $. The graph of $ \nabla F $ defines a special
Lagrangian submanifold $ M $ in $ \mathbb C^n = \mathbb R^n \times \mathbb R^n $.
In other words,
$ F $ satisfies equation (\ref{eq1.1}). If there is a constant
$ \beta < \infty $ that satisfies (\ref{eq1.3}) and (\ref{eq1.4}), then $ F $ is a quadratic
polynomial and $ M $ is flat. 
\end{theo}

The celebrated theorem of Bernstein says that the
only entire minimal graphs in Euclidean $ 3 $-space are
planes. This result has been partially ge\-ne\-ra\-lized to 
higher codimension. If $ f : \mathbb R^k \to \mathbb R^n $ is an entire
solution of the minimal surface system with bounded
gradient, then $ f $ is linear for $ k = 2 $ by a theorem of
Osserman-Chern and for $ k = 3 $ by a result of Fischer-Colbrie \cite{FC}. 
For larger $ k $, however, there exist counterexamples of
Lawson-Osserman. By way of contrast, our
Theorem \ref{theo1.3} shows that minimal Lagrangian
graphs with bounded gradient are always planar. 

\section{Geometry of Lagrangian Grassmannian manifolds}

Let $ \mathbb R^{m+n} $ be an $ (m+n) $-dimensional Euclidean space. The set of all oriented
$ n $-subspaces (called $ n $-planes) constitutes the Grassmannian manifold $ G_{n,m}  $, which
is the irreducible symmetric space $ SO (m+n) / SO (m) \times SO (n) $.

Let $ \{ e_{\alpha}, e_{n+i} \}$ be a local orthonormal frame field in $ \mathbb R^{m+n} $, where 
$ i, j, ... = 1, ..., m; \alpha, \beta, ... = 1, ..., n; a, b, ... = 1, ..., m + n $ (say, 
$ n \leq m$). Let 
$ \{ \omega_{\alpha} , \omega_{n + i} \} $ be its dual frame field so that the Euclidean metric is
$$
 g = \sum_{\alpha} \omega^2_{\alpha} + \sum_i \omega^2_{n+i} .
$$
The Levi-Civita connection forms $ \omega_{ab} $ of $ \mathbb R^{m+n} $ are uniquely determined by
the equation
\begin{align}
\label{eq2.1}
 d \omega_a & = \omega_{ab} \wedge \omega_b , \\ \nonumber
 \omega_{ab}       & + \omega_{ba} = 0 .
\end{align}
The canonical Riemannian metric on $ G_{n,m} $ can be defined by
\begin{equation}
\label{eq2.2}
 ds^2 = \sum_{\alpha, i} \omega^2_{\alpha \, n + i} .
\end{equation}
>From (\ref{eq2.1}) and (\ref{eq2.2}) it is easily seen that the curvature tensor of $ G_{n,m} $ is
\begin{equation}
\label{eq2.3}
 \begin{split}
	R_{\alpha i \, \beta j \, \gamma k \, \delta l} = \delta_{\alpha \beta} \delta_{\gamma \delta}
	 \delta_{ik} \delta_{jl} + \delta_{\alpha \gamma} \delta_{\beta  \delta } \delta_{ij} \delta_{kl} \\
	 - \delta_{\alpha \beta} \delta_{\gamma \delta} \delta_{il} \delta_{kj} - \delta_{\alpha \delta}
 		\delta_{\beta \gamma} \delta_{ij} \delta_{kl}
 \end{split}
\end{equation}
in a local orthonormal frame field $ \{ e_{\alpha i}  \} $, which is dual to $ \{ \omega_{\alpha \, n+1} \}$.\\

Let $ P_0 $ be an oriented $ n $-plane in $ \mathbb R^{m+n} $. We represent it by $ n $ vectors $ e_{\alpha} $,
which
are complemented by $ m $ vectors $ e_{n+i} $, such that $ \{ e_{\alpha}, e_{n+i} \} $ form an orthonormal
base of $ \mathbb R^{m+n} $. Then we can span the $ n $-planes $ P $ in a neighborhood $ \mathbb U $ of $ P_0 $
by $ n $ vectors $ f_{\alpha} $:
$$
 f_{\alpha} = e_{\alpha} + z_{\alpha z} e_{e+i}, 
$$
where $ (z_{\alpha i}) $ are the local coordinates of $ P $ in $ \mathbb U $. The metric (\ref{eq2.2}) on
$ G_{n,m} $ in those local coordinates can be described as
\begin{equation}
\label{eq2.4}
 ds^2 = tr (( I_n + ZZ^T )^{-1} dZ (I_m + Z^TZ)^{-1} dZ^T ) 
\end{equation}
where $ Z = (z_{\alpha i}) $ is an $ (n \times m) $-matrix and $ I_n $ (res. $ I_m $) denotes the
$ (n\times n) $-identity (res. $ m \times m $) matrix. 

Now we consider the case $ \mathbb R^{2n} = \mathbb C^n $ which has the
usual complex structure $ J $. For any $ u = (x, y) = (x_1, ... x_n; y_1, ... y_n) $ in 
$ \mathbb R^{2n} $
$$
 J u = ( - y_1, ..., -y_n ; x_1, ..., x_n) .
$$
For an $ n $-plane $ \zeta \subset \mathbb R^{2n} $ if any $ u \in \zeta $ satisfies
$$
 \langle u, J u \rangle = 0 ,
$$
then $ \zeta $ is called a Lagrangian plane. 
The Lagrangian planes yield the Lagrangian Grassmannian
manifold. It is the symmetric space $ \du{U(n)}{SO(n)} $.\\

For the Grassmannian $ G_{n,n} $ its local coordinates
are described by $ (n \times n) $-matrices. For any $ A \in LG_n \subset G_{n,n} $,
let $ u = (x, x A) $ and $ \tilde u = (\tilde x, \tilde x A) $ be two vectors in $ A $. 
By definition $ \langle u, J \tilde u \rangle = 0 $ and so we have
$$ 
 A = A^T
$$
>From (\ref{eq2.4}) it is easy to see that the transpose is
an isometry of $ G_{n,n} $.
Hence the fix point set $ LG_n $ is a totally
geodesic submanifold of $ G_{n,n} $.
By the Gauss equation the Riemannian curvature
tensor of $ LG_n $ is also defined by (\ref{eq2.3}).\\

Let $ \dot \gamma = x_{\alpha i} e_{\alpha i} $ be a unit tangent vector at $ P_0 $, where
$ \{ e_{\alpha i} \} $ is a local orthonormal frame field. By an 
action of $ SO (n) $
$$
 x_{\alpha i} = \lambda_{\alpha} \delta_{\alpha i} ,
$$
there $ \sum\limits_{\alpha} \lambda^2_{\alpha} = 1 $. In our previous paper \cite{JX} we have
computed the eigenvalues of the Hessian of the distance
function from a fixed print $ P_0 $ at the direction
$ \dot \gamma = (x_{\alpha i}) = (\lambda_{\alpha} \delta_{\alpha i}) $.
Considering the present situation, 
when $ m = n $ and the eigenvectors are symmetric matrices,
the eigenvalues are as follows:
\begin{align}
\label{eq2.5}
	(\lambda_{\alpha} - \lambda_{\beta}) \cot (\lambda_{\alpha} - \lambda_{\beta} ) r
 	& \quad \hbox{with multiplicity } \; 1 \\ \nonumber
	\frac{1}{r} 
	& \quad \hbox{with multiplicity } \; n - 1
\end{align}
where $ r $ is the distance from $ P_0 $, for sufficiently small $ r $.\\

The geodesic from $ P_0 $ at $ (x_{\alpha i}) = (\lambda_{\alpha} \delta_{\alpha i}) $ in the local
coordinates neighborhood $ U $ is (see \cite{W})

\begin{equation}
\label{eq2.6}
(z_{\alpha i} (t) ) = \begin{pmatrix} \tan (\lambda_1 t) &  & 0 \\
                                   & \ddots &  \\ 
                                 0 &  & \tan (\lambda_n t) 
                               \end{pmatrix}
\end{equation}
where $ t $ is the arc length parameter and $ 0 \leq t < \frac{\pi}{2 | \lambda_n |} $ with
$ | \lambda_n | = \max ( | \lambda_1 |, ... , | \lambda_n |  ) $.

\section{Gauss map}

Let $ M $ be an $ n $-dimensional oriented submanifold in $ \mathbb R^{m+n} $. Choose an 
orthonormal frame field $ \{e_1, ..., e_{m+n} \} $ in $ \mathbb R^{m+n} $ such that the $ e'_{\alpha} s $
are tangent to $ M $. Let $ \{ \omega_1, ..., \omega_{m+n}  \} $ be its coframe field. 
Then, the structure equations of $ \mathbb R^{m+n} $ along $ M $ are as 
follows. 
$$
\omega_{n+i} = 0,
$$
$$
d \omega_{\alpha} = \omega_{\alpha\beta}\wedge \omega_{\beta}, \, \omega_{\alpha\beta} + \omega_{\beta\alpha} = 0, 
$$
$$
\omega_{n+i \, \alpha} = h_{i \alpha \beta} \omega_{\beta}, 
$$
$$
d \omega_{ij} = \omega_{i \kappa} \wedge \omega_{\kappa j} + \omega_{i \alpha} \wedge \omega_{\alpha j}, 
$$
where the $ h_{i \alpha \beta} $, the coefficients of the second fundamental form of $M  $ in $ \mathbb R^{m+n} $, 
are symmetric in $ \alpha $ and $ \beta $.
Let $ 0 $ be the origin of $ \mathbb R^{m+n} $. Let $ SO (m+n)$ be the manifold consisting of
all the orthonormal frames $ (0; e_{\alpha}, e_{n+i}) $. Let $ P = \{ (x; e_1, ..., e_n ); x \in
M, e_{\alpha} \in T_x M \} $ be the principal bundle of orthonormal tangent frames over $ M, Q =
\{ ( x; e_{n+1}, ..., e_{m+n}  ) ; x \in M, e_{n+i} \in N_x M \} $ be the principal bundle of orthonormal
normal frames over $ M$, then $ \bar\pi : P \otimes Q \to M $ is the projection with fiber
$ SO (m) \times SO (n), i: P \otimes Q \hookrightarrow SO (m+n)  $ is the natural inclusion. \\

We define the generalized Gauss map $ \gamma : M \to {\bf G}_{n,m} $ by
$$
 \gamma (x) = T_x M \in {\bf G}_{n,m}
$$
via the parallel translation in $ \mathbb R^{m+n} $ for $ \forall x \in M $. Thus, the following commutative 
diagram holds
\begin{equation*}
\begin{CD}
 P \otimes Q @>i>> SO (m+n)  \\
 @V{\bar\pi}VV     @VV{\pi}V \\
 M  @>{\gamma}>>  {\bf G}_{n,m} 
\end{CD}
\end{equation*}
Using the above diagram, we have
\begin{equation}
\label{eq3.1}
 \gamma^{\ast} \omega_{n+i \, \alpha} = h_{i \alpha \beta } \omega_{\beta}.
\end{equation}
Now, we assume that $ M $ is a Lagrangian submanifold
in $ \mathbb R^{2n} $. The image of the Gauss map $ \gamma : M \to {\bf G}_{n,m} $ then lies in its 
Lagrangian Grassmannian $ LG_n $.
We then have $ e_{\alpha} \in TM $ and $ J e_{\alpha} \in NM $.\\

Furthermore, 
$$
 h_{i \alpha \beta} = 
 \langle \nabla_{e_{\alpha}} e_{\beta}, J e_i  \rangle =
 - \langle \nabla_{e_{\alpha}} J e_i, e_{\beta}  \rangle =
 \langle \nabla_{e_{\alpha}} e_i, J e_{\beta}  \rangle =
 h_{\beta \alpha i} .
$$
Thus, the $ h_{i \alpha \beta} $ are symmetric in all their indices.
(\ref{eq3.1}) can also be written as its dual form
\begin{equation}
\label{eq3.2}
 \gamma_{\ast} e_{\beta} = h_{i \alpha \beta} e_{\alpha i} .
\end{equation}
For each $ e_{\beta} , (h_{i \alpha \beta} )$ is a symmetric matrix. 

\section{Minimal Legendrian submanifolds in the sphere and minimal Lagrangian cones}

\indent In the sphere $ S^{2n-1} \hookrightarrow \mathbb R^{2n }$ there is a standard contact
structure. Let $ X $ be the position vector field
of the sphere and $ \eta $ be the dual form of $ J X $ in $ S^{2n-1} $, where
$ J $ is the complex structure of $ \mathbb C^n = \mathbb R^{2n} $. It 
is easily seen that
\begin{equation}
\label{eq4.1}
d \eta = 2 \omega ,
\end{equation}
where $ \omega $ is the K\"ahler form of $\mathbb C^n  $. Therefore,
\begin{equation}
\label{eq4.2}
\eta \wedge (d \eta)^{n-1} \neq 0
\end{equation}
everywhere and $ \eta $ is a constant form in $ S^{2n-1} $.
The maximal dimensional integral submanifolds of 
the distribution 
\begin{equation}
\label{eq4.3}
 \eta = 0
\end{equation}
are $ (n-1) $-dimensional and are called Legendrian 
submanifolds in $ S^{2n-1} $.\\

Now, let us consider the cone $ CM  $ over $ M $. 
$CM$ is the image under the map $ M \times [ 0, \infty ) $ into $ \mathbb R^{2n} $
defined by $ (x, t)  \to t x $, where $ x \in M, t \in [ 0, \infty )$.
$ CM $ has a singularity at $ t = 0 $. The associated
truncated cone $ CM_{\varepsilon} $ is the image of $  M \times [ \varepsilon, \infty )$
under the same map, where $ \varepsilon $ is any positive number.\\

We have (see \cite{S})
\begin{prop1}
\label{prop4.1}
$ CM_{\varepsilon} $ is minimal submanifold
in $ \mathbb R^{2n} $if and only if $ M $ is a minimal submanifold in 
$ S^{2n-1} $. 
\end{prop1}
For a fixed point $ x \in M$ choose a local orthonormal 
frame field $ \{ e_s \} $ ($s = 1, ..., n-1$) near $ x $ in $ M $ with 
$ \nabla_{e_s} e_t |_x = 0$.\\

By parallel translating along rays from the origin,
we obtain a local vector field $ E_s $ in $ CM $. Obviously, 
$ E_s = \frac{1}{r} e_s $, where $ r $ is the distance from the origin. 
Thus, $ \{ E_s, \tau  \} $ is a frame field in $CM  $, where $ \tau = \frac{\partial}{\partial r} $
is the unit tangent vector along rays.
Obviously $ \nabla_{\tau} \tau = 0 $.

In the case of $ M $ being Legendrian
$$
 \eta (e_s) = 0 ,
$$
and 
$$
 d \eta (e_s, e_t) = ( \nabla_{e_s} \eta )e_t - (\nabla_{e_t} \eta ) e_s 
$$
$$
 = \nabla_{e_s} \eta (e_t) - \nabla_{e_t} \eta (e_s) - \eta ( [ e_s, e_t ] ) = 0 .
$$
>From (\ref{eq4.1}) it follows that 
\begin{equation}
\label{eq4.4}
 \omega (E_s, E_t) = \frac{1}{r^2} \omega ( e_s, e_t ) = 0.
\end{equation}
Obviously
\begin{equation}
\label{eq4.5}
 \omega (E_s, \tau) = \langle E_s, J \tau  \rangle = \frac{1}{r} \eta (e_s) = 0
\end{equation}
(\ref{eq4.4}) and (\ref{eq4.5}) mean that $ CM $ is a Lagrangian submanifold
in $ \mathbb R^{2n} $ if and only if $ M $ is a Legendrian submanifold
in $ S^{2n-1} $. \\

Now, let us compute the coefficients of the second
fundamental form of $ CM $ in $ \mathbb R^{2n} $.\\

We have a local orthonormal frame 
field $ \{ E_s, \tau, J E_s, J \tau  \} $ in $ \mathbb R^{2n} $ along $ CM $, where
$ \{ E_s, \tau \} $ is a local orthonormal frame field in $ CM $. 

Note
$$
 \nabla_{E_s} \tau = \nabla_{E_s} \frac{X}{r} = \frac{1}{r} E_s ,
$$
where $ X $ denotes the position vector of the concernd point.

Then
$$
 \langle \nabla_{E_s} E_t, \tau  \rangle = - \langle E_s, \nabla_{E_t} \tau \rangle  = - \frac{1}{r}
 \delta_{st} ,
$$
and
\begin{align*}
\frac{d}{d r} \langle \nabla_{E_s} E_t , E_u \rangle 
 & = \langle  \nabla_{\tau } \nabla_{E_s} E_t , E_u \rangle \\
 & = \langle \nabla_{E_s } \nabla_{\tau } E_t, E_u  \rangle + \langle \nabla_{[\tau, E_s] } E_t, E_u  \rangle \\
 & = - \frac{1}{r } \langle \nabla_{E_s } E_t, E_u  \rangle , \\
\frac{d}{d r} \langle \nabla_{E_s} E_t , J E_u \rangle 
 & = - \frac{1}{r} \langle \nabla_{E_s} E_t , J E_u \rangle .
\end{align*}
Integrating them gives
$$
 \langle \nabla_{E_s} E_t , E_u \rangle = \frac{C_{ ust}}{r}
$$
and
$$
 \langle \nabla_{E_s} E_t , J E_u \rangle = \frac{ D_{ ust}}{ r} ,
$$
where $ C_{ ust} , D_{ ust} $ are constants along the ray. They can
be determined by the conditions at $ r = 1 $ as follows
$$
 C_{ ust} = 0 \quad , \quad D_{ ust} = h_{ust} ,
$$
where $ h_{ust} $ are the coefficients of the second fundamental form of $ M $ in $ S^{2n-1} $
in the $ J e_u $ directions. We also have 
$$
\langle \nabla_{E_s} E_t, J \tau   \rangle  = 
 - \langle E_t, \nabla_{E_s} J \tau \rangle =
 - \frac{1}{r} \langle E_t, J E_s  \rangle  = 0.
$$
Thus, we obtain the coefficients of the second fundamental form
$ CM $ in $ \mathbb R^{2n} $ as follows. 
In the $ J E_u $ directions 
\begin{equation}
\label{eq4.6}
 B_{uij} = \begin{pmatrix} \frac{h_{uij}}{r} & 0 \\
                         0 & 0                    \end{pmatrix}
\end{equation}
and in the $ J \tau $ direction
\begin{equation}
\label{eq4.7}
 B_{nij} = 0
\end{equation}
>From (\ref{eq3.2}), (\ref{eq4.5}) and (\ref{eq4.6}) we know that the
Gauss map of the cone $ CM $ has rank $ n-1 $ at most. 
We summarize the results of this section as 
\begin{prop1}
\label{prop4.2}
Let $ M $ be an $ (n-1) $-dimensional submanifold in $ S^{2n-1} $. It is minimal and Legendrian if
and only if the cone $ CM $ over $ M  $ is
a minimal Lagrangian submanifold in
$ \mathbb R^{2n} $. Furthermore, the Gauss map $ \gamma : CM \to LG_n $
has rank $ n-1 $ at most.
\end{prop1}

\section{Harmonic Maps}

Let $ (M, g) $ and $ (N, h) $ be Riemannian manifolds with metric tensors $ g $ and $ h $,
respectively. Harmonic maps are described as critical points of the following energy
functional 
\begin{equation}
\label{eq5.1}
 E (f) = \frac{1}{2} \int\limits_M e (f) \ast 1 ,
\end{equation}
where $ e (f)  $ stands for the energy density. The Euler-Lagrange equation of the 
energy functional is 
\begin{equation}
\label{eq5.2}
 \tau (f) = 0 ,
\end{equation}
where $ \tau (f)  $ is the tension field. In local coordinates
\begin{align}
\label{eq5.3}
e (f) 
& = g^{if} \frac{\partial f^{\beta}}{\partial x^i} \frac{\partial f^{\gamma}}{\partial x^j} h_{\beta \gamma} , \\
\label{eq5.4}
\tau (f) 
& = ( \Delta_M f^{\alpha} + g^{ij} \Gamma^{\alpha}_{\beta \gamma } 
      \frac{\partial f^{\beta}}{\partial x^i} \frac{\partial f^{\gamma}}{\partial x^j} )  
      \frac{\partial}{\partial y^{\alpha}} ,
\end{align}
where $ \Gamma^{\alpha}_{\beta \gamma} $ denotes the Christoffel symbols of the target manifold $ N $. 
For more details on harmonic maps consult \cite{EL}.

A Riemannian manifold $ M $ is said to be simple, if it can be described by coordinates $ x $
on $ \mathbb R^n $ with a metric
\begin{equation}
\label{eq5.5}
 ds^2 = g_{ij} d x^i dx^j ,
\end{equation}
for which there exist positive numbers $ \lambda $ and $ \mu $ such that 
\begin{equation}
\label{eq5.6}
 \lambda | \xi |^2 \leq g_{ij } \xi^i \xi^j \leq \mu | \xi |^2
\end{equation}
for all $ x $ and $ \xi $ in $ \mathbb R^n $. In other words, $ M $ is topologically $ \mathbb R^n $
with a metric for
which the associated Laplace operator is uniformly elliptic on $ \mathbb R^n $.

Hildebrandt-Jost-Widman proved a Liouville-type theorem for harmonic maps in \cite{HJW}: 
\begin{theo1}
\label{theo5.1}
Let $ f $ be a harmonic map from a simple or compact Riemannian
manifold $ M $ into a complete Riemannian manifold $ N $, the sectional curvature of
which is bounded above by a constand $ \kappa \geq 0 $. Denote by $ B_R (Q)  $ a geodesic ball in
$ N $ with radius $ R < \frac{\pi}{2 \sqrt{\kappa}} $ which does not meet the cut locus of its center $ Q $.
Assume also that the range $ f (M) $ of the map $ f $ is contained in $ B_R (Q) $. Then $ f $ is a constand
map. 
\end{theo1}
\begin{rem}
In the case where $ B_R (Q) $ is replaced by another geodesically convex 
neighborhood, the iteration technique in \cite{HJW} is still applicable and the result
remains true (for example, a general version of that iteration technique that
directly applies here has been given in \cite{GJ}).
\end{rem}
By using the composition formula for the tension field, one easily verifies that
the composition of a harmonic map $ f : M \to N $ with a convex function $ \phi: f (M) \to \mathbb R $
is a subharmonic function on $ M $. The maximum principle then implies
\begin{prop1}
\label{prop5.2}
Let $ M $ be a compact manifold without boundary, $ f : M \to N $
a harmonic map with $ f (M) \subset V \subset N $. Assume that there exists a strictly convex
function on $ V $. Then $ f $ is a constant map.
\end{prop1}
Let $ M \to \mathbb R^{m+n} $ be an $ n $-dimensional oriented submanifold in
Euclidean space. We have the relation between the property of the submanifold
and the harmonicity of its Gauss map in \cite{RV}.
\begin{theo1}
\label{theo5.3}
Let $ M  $ be a submanifold in $ \mathbb R^{m+n} $. Then the mean curvature vector
of $ M $ is parallel if and only if its Gauss map is a harmonic map. 
\end{theo1}

Let $ M \to S^{m+n} \hookrightarrow \mathbb R^{m+n-1}$ be an $ m $-dimensional submanifold in the sphere. 
For any $ x \in M $, by parallel translation in $ \mathbb R^{m+n+1} $, the normal space $ N_x M $ of
$M$ in $ S^{m+n} $ is moved to the origin of $ \mathbb R^{m+n+1} $. We then obtain an $ n $-subspace
in $ \mathbb R^{m+n+1}  $. Thus, the so-called normal Gauss map $ \gamma : M \to G_{n,m+1} $ has been
defined. There is a natural isometry $ \eta $ between $ G_{n, m+1} $ and $ G_{m+1,n} $ which maps
any $ n $-subspace into its orthogonal complementary $ (m+1) $-subspace. The map
$ \eta^{\ast} = \eta \circ \gamma$ maps any point $ x \in M $ into an $ (m+1) $-subspace spanned by 
$ T_x M $ and
the position vector of $ x $.
>From Theorem \ref{theo5.3} and Proposition \ref{prop4.1} it follows that
\begin{prop1}
\label{prop5.4}
$ M $ is a minimal $ m $-dimensional submanifold in the sphere
$ S^{m+n} $ if and only if its normal Gauss map $ \gamma : M \to G_{n, m+1} $ is a harmonic map.
\end{prop1}
\section{Proofs of the theorems}
Proof of Theorem 1 \\

Since $ M $ is a graph in $ \mathbb R^{2n} $ defined by $ \nabla F $, the
induced metric $ g $ on $ M $ is 
$$
d s^2 = g_{\alpha \beta } dx^{\alpha} dx^{\beta},
$$
where 
$$
 g_{\alpha \beta} = \delta_{\alpha \beta} + \frac{\partial^2 F}{\partial x^{\alpha } \partial x^{\gamma}}
 \frac{\partial^2 F}{\partial x^{\beta} \partial x^{\gamma}}. 
$$
It is obvious that the eigenvalues of the matrix $ (g_{\alpha \beta}) $
at each point are $ \geq 1 $. The condition (\ref{eq1.3}) implies
that the eigenvalues of the matrix $ (g_{\alpha \beta}) $ are $ \leq \beta^2 $.

The condition (\ref{eq5.6}) is satisfied and $ M $ is a simple
Riemannian manifold.\\

Let $ \{ e_{\alpha}, e_{n+\beta} \} $ be the standard orthonormal base
of $ \mathbb R^{2n} $. Choose $ P_0 $ as an $ n $-plane spanned by $ e_1 \wedge ... \wedge e_n $.
At each point in $ M $ its image $ n $-plane $ P $ under the
Gauss map is spanned by
$$
 f_{\alpha}
 = e_{\alpha} + \frac{\partial^2 F}{\partial x^{\alpha } \partial x^{\beta}} e_{n+\beta},
$$
which lies in the Lagrangian Grassmannian manifold $ LG_n $.

Suppose the eigenvalues of Hess $(F) $ at each point $  x$ are
$ \mu_{\alpha} (x) $ which are positive by the convexity of the function $F$.
The condition (\ref{eq1.3}) means
$$
\prod\limits_{\alpha} (1 + \mu^2_{\alpha} ) \leq \beta .
$$
Hence, 
\begin{equation}
\label{eq6.1}
\mu_{\alpha} \leq \sqrt{\beta^2 - 1}
\end{equation}
Define in the normal polar coordinates of $ P_0 $ in $ LG(n) $
$$
 \tilde B_{LG} (P_0) = \{ (X, t) ; X= (\lambda_{\alpha}\delta_{\alpha i} ) ,
 \lambda_{\alpha } \geq 0; 
 0 \leq t \leq t_x = \tan \sqrt[-1]{\beta^2 -1} \} .
$$
Two points $ P_0 $ and $ P  $ can be joined by a unique geodesic
$ P(t) $ spanned by 
$$
 \tilde f_{\alpha} (t) = e_{\alpha} + z_{\alpha \beta }(t) e_{n+\beta}  ,
$$
where
$$
z_{\alpha \beta } (t) = \begin{pmatrix} \tan (\lambda_1 t) & & 0 \\
                                                              & \ddots & \\
                                           0                  & & \tan ( \lambda_n t)
                           \end{pmatrix} .
$$
Therefore, the image under the Gauss map $\gamma  $ of $ M $ lies in
$ \tilde B_{LG} (P_0) $. \\

On the other hand, from (\ref{eq2.5}) we see that when $ \lambda_{\alpha} \geq 0$
the square of the distance function $ r^2 $ from $ P_0 $ is
a strictly convex smooth function in $ \tilde B_{LG} (P_0)$.
Furthermore, it is a geodesically convex set. 

Now, we have the Gauss map $ \gamma : M \to \tilde B_{LG} (P_0) \subset LG_n$ 
which is harmonic by Thoerem \ref{theo5.3}.
Hence the conclusion follows by using Theorem \ref{theo5.1}.
\begin{rem}
If the graph of $ \nabla F $ is a submanifold
with parallel mean curvature instead of a minimal submanifold
the Theorem remains true as well. 
\end{rem}

Proof of Theorem 2 \\

>From Proposition \ref{prop5.4} we know that the normal Gauss map
$ \gamma : M \to G_{n,n} $ is harmonic. Let $\eta $ be the 
isometry in $ G_{n,n} $ which maps any $ n $-plane
into its orthogonal complementary $ n $-plane. Hence 
$\eta \circ \gamma  $ is also harmonic. On the other hand, from
the discussion in \S~4 it follows that
\begin{equation}
\label{eq6.2}
\eta \circ \gamma (M) = \gamma' (CM_{\varepsilon}),
\end{equation}
where $ \gamma' : CM_{\varepsilon} \to G_{n,n}$ is the Gauss map. Now
let $ \eta (P_0 ) $ be spanned by $ n $ vectors $ e_{\alpha} $, which are complemented
by $ n $ vectors $ e_{n+i} $. The condition (\ref{eq1.5})
ensures that for all normal $ n $-planes $ P $ of $ M $ in $S^{2n-1}, \eta (P)  $ lies
in the coordinate neighborhood $ U $ of $ \eta (P_0) $
and $ \eta (P) $ is spanned by $n $ vectors $ f_{\alpha} $:
$$
 f_{\alpha} = e_{\alpha} + z_{\alpha i} e_{n+i} ,
$$
where $ (z_{\alpha i}  )$ are local coordinates of $ \eta (P) $ in $U $. Noting
(\ref{eq6.2}), $ \eta (P) $ lies in $ LG_n $ and $ (z_{\alpha i}) $ is a symmetric
matrix. the geodesic
from $ \eta (P_0) $  at $ (X_{\alpha i}) = (\lambda_{\alpha } \delta_{\alpha i})$ in $U $ is
described by (\ref{eq2.6}).

Noting Proposition \ref{prop4.2} there exists $ \alpha_0$ such that $ \lambda_{\alpha_0} = 0$.
Then by actions of $ SO (n) $ we can achieve that at most
one of $ \lambda_{\alpha}'s $ is negative. Then by
one more action of $ SO (n) $ all the $ \lambda_{\alpha}'s $ are
nonnegative. 

Take any point $ \eta (p) $ in $ \eta \circ \gamma (M) = \gamma' (CM_{\varepsilon}) $. Draw a 
geodesic $ P (t) $ from $\eta (P_0)  $ to $ \eta (P) $. Let $ P (t) $ be spanned 
by
$$
 f_{\alpha} = e_{\alpha} + z_{\alpha i} e_{n+i}
$$
where $ z_{\alpha i} $ is defined by (\ref{eq2.6}). Let 
$$
 \tilde f_1 = \cos ( \lambda_1 t) f_1, ... , \tilde f_n = \cos ( \lambda_n t) f_n .
$$
Those $ \tilde f_1, ..., \tilde f_n $ are orthonormal. Therefore,
$$
 \langle \eta (P_0), P (t) \rangle = \prod\limits^n_{\alpha = 1} \cos (\lambda_{\alpha} t) ,
$$
where $ \lambda_{\alpha} \geq 0 $ and $ \Sigma \lambda^2_{\alpha} = 1 $. From (\ref{eq1.5}) it follows
that
\begin{equation}
\label{eq6.3}
 t \leq \frac{\cos^{-1} \delta}{\max\limits_{\alpha} (\lambda_{\alpha})}.
\end{equation}
Define in the normal polar coordinates around $ \eta (P_0) $
\begin{equation}
\label{eq6.4}
 \tilde B_{LG} ( \eta (P_0) ) = 
 \{ (X, t) ; X = (\lambda_{\alpha } \delta_{\alpha i}), 
 0 \leq t \leq t_X = 
 \frac{\cos^{-1}\delta}{\max\limits_{\alpha} (\lambda_{\alpha})}
 \} .
\end{equation}
>From (\ref{eq2.6}) we see that $ \tilde B_{LG} ( \eta (P_0)) $ lies inside the cut locus of $ \eta (P_0) $. 
We also know from (\ref{eq2.5}) that the square of the distance
function $ r^2 $ from $ \eta (P_0) $ is a strictly convex smooth
function in $ \tilde B_{LG} (\eta (P_0)) $.  By a similar argument as for
$ B_G (P_0) $ in our previous paper \cite{JX} it can be shown that $\tilde B_{LG} (P_0)  $ is a geodesically
convex set.\\

We thus have a harmonic map $ \eta \circ \gamma $ from
$ M $ into a geodesically convex set $ \tilde B_{LG} (\eta (P_0)) $. By 
using Theorem \ref{theo5.1} we conclude that $ \eta \circ \gamma $ is a
constant map, and then so is the map $ \gamma $. This completes
the proof. \\

Proof of Theorem 3 \\

Let us consider the tangent cone of $ M $ at $ \infty $ as Fleming in \cite{Fl}.
Take the intersection of $ M $ with the ball of radius
$ t $ and contract by $ \frac{1}{t} $ to get a family of
minimal submanifolds in the unit ball with 
submanifolds of $ S^{2n-1} $ as boundaries. More precisely, we define a sequence
$$
 F^t = \frac{1}{t^2} F(tx) .
$$
For each $ t $
$$
 \frac{\partial^2 F^t}{\partial x^{\alpha} \partial x^{\beta}} =
 \frac{\partial^2 F}{\partial u^{\alpha} \partial u^{\beta}} ,
$$
where $ u^{\alpha} = tx^{\alpha}$. It turns out $ F^t  $ satisfies the same 
conditions as $ F $. Moreover, there 
is a subsequence $ t_j \to \infty $ such that
$$
 \lim\limits_{t_j \to \infty} F^t (x) = \tilde F (x).
$$
$ \tilde F $ satisfies (\ref{eq1.1}), (\ref{eq1.2}) and (\ref{eq1.3}) and the graph $ \nabla \tilde F $ is
a special Lagrangian cone $ C \tilde M $ whose link 
is a compact minimal Legendrian submanifold $ \tilde M $.

Let $ \{ e_{\alpha}, e_{n+\beta}  \} $ be the standard orthonormal base 
of $ \mathbb R^{2n} $. Choose $ P_0 $ as an $ n $-plane spanned by $ e_1 \wedge ... \wedge e_n $.
At each point of $ C \tilde M $ its image $ n $-plane $ P $ under
the Gauss map is spanned by
$$
 f_{\alpha} = e_{\alpha } + \frac{\partial^2 \tilde F}{\partial x^{\alpha} \partial x^{\beta }} e_{n+\beta}
$$
It follows that 
$$
 | f_1 \wedge ... \wedge f_n  |^2 = \det \left( \delta_{\alpha \beta} + 
 \frac{\partial^2 \tilde F}{\partial x^{\alpha} \partial x^{\gamma }} 
 \frac{\partial^2 \tilde F}{\partial x^{\beta} \partial x^{\gamma }} \right)
$$
and 
$$
 \Delta_f = | f_1 \wedge ... \wedge f_n  |.
$$
The $ n $-plane $ P $ is also spanned by
$$
 P_{\alpha} = \Delta^{-\frac{1}{n}}_f f_{\alpha}, 
$$
moreover,
$$
 | p_1 \wedge ... \wedge p_n  | = 1 .
$$ 
We then have 
\begin{align*}
\langle P, P_0 \rangle & = \det ( \langle e_{\alpha}, P_{\beta} \rangle  ) \\
                        & = \Delta^{-1}_f \geq \beta^{-1} .
\end{align*}
Let $ \eta $ be the isometry in $ G_{n,n} $ that 
maps any $ n $-plane into its orthogonal complementary 
$ n $-plane. We thus have
$$
 \langle \eta P , \eta P_0 \rangle \geq \beta^{-1} .
$$  
By the discussion in \S~4 we know that $ \eta P $ is
just the normal $ n $-plane of $ \tilde M $ in $ S^{2n-1} $. Then Theorem 2
tells
us that $ \tilde M$ is a totally geodesic sphere $ S^{n-1} $ in $ S^{2n-1} $ and therefore,
$ C \tilde M $ is an $ n $-plane in $ \mathbb R^{2n} $. Allard's result \cite{A}
then implies that the original special Lagrangian submanifold
$ M $ is an affine $ n $-plane and $ F $ is a quadratic polynomial.


\newpage
\begin{thebibliography}{99} 
\bibitem{A}	Allard, W., 
		On the first variation of a varifold, 
		Ann. Math. 95 (1972), 417-491.
\bibitem{B}	Barbosa, J.L.M., 
		An extrinsic rigidity theorem for minimal immersions from $ S^2$ into $ S^n $,
		J. Differential Geometry 14(3) (1980), 355-368.
\bibitem{C}	Calabi, E., 
		Improper affine hyperspheres of convex type and genera\-li\-za\-tion of a theorem by K. J\"orgens,
		Mich. Math. J. 5 (1958), 105-126.
\bibitem{CDVV}	Chen, B.Y., Dillen, F., Verstraelen, L. and Vrancken, L.,
		 An exotic real minimal immersion of $ S^3 $ in $ CP^3 $ and its characterisation,
		Proc. Roal Soc. Ed. 126(A) (1996), 153-165.
\bibitem{EL}	Eells, J. and Lemaire, L., 
		Another report on Harmonic maps, 
		Bull. London Math. Soc. 20(5) (1988), 385-524. 
\bibitem{Fl}	Fleming, W.H., 
		On the oriented Plateau problem,
		Circolo Mat. Palermo II (1962), 1-22. 
\bibitem{Fu}	Fu, Lei,
		An analogue of Bernstein's theorem, 
		Houston J. Math. 24(3) (1998), 415-419. 
\bibitem{FC}	Fischer-Colbrie, D., 
		Some rigidity theorems for minimal submanifolds of the sphere,
		Acta math. 145 (1980), 29-46.  
\bibitem{GJ}	Gulliver, R. and Jost, J., 
		Harmonic maps which solve a free boundary problem,
		J. Reine Angew. Math. 381 (1987), 61-89.  
\bibitem{H} 	Haskins, Mark,
		Special Lagrangian cones,
		math. DG/0005164.
\bibitem{HL}	Harvey, R. and Lawson, H.B.,
		Calibrated geometry,
		Acta Math. 148 (1982), 47-157.  
\bibitem{HJW}	Hildebrandt, S., Jost, J. and Widman, K.O., 
		Harmonic mappings and minimal submanifolds,
		Invent. math. 62 (1980), 269-298.
\bibitem{JX}	Jost, J. and Xin, Y.L., 
		Bernstein type theorems for higher codimension, 
		Calc. Var. 9 (1999), 277-296.  
\bibitem{LO}	Lawson, B. and Osserman, R.,
		Non-existence, non-uniqueness and irregularity of solutions to the minimal surface system,
		Acta math. 139 (1977), 1-17.  
\bibitem{P}	Pogorelov, A.V., 
		On the improper affine hypersurfaces, 
		Geom. Dedicata 1 (1972), 33-46.
\bibitem{RV}	Ruh, E.A. and Vilms, J., 
		The tension field of the Gauss map,
		Trans. A.M.S. 149 (1970), 569-573.  
\bibitem{S}	Simons, J., 
		Minimal varieties in Riemannian manifolds,
		Ann. Math. 88 (1968), 62-105.  
\bibitem{SYZ}	Strominger, A., Yau, S.T. and Zaslow, E., 
		Mirror symmetry is T-duality, 
		Nucl. Phys. B479 (1996), 243-259.
\bibitem{TW}	Trudinger, Neil S. and Wang, Xu-Jia,
		The Bernstein problem for affine maximal hypersurfaces,
		Invent. math. 140 (2000), 399-422.  
\bibitem{W}	Wong, Yung-Chow, Differential geometry of Grassmann manifolds, 
		Proc. N.A.S. 57 (1967), 589-594.  
\end{thebibliography}
\end{document}